# Proof of aperiodicity of hat tile using the golden ratio

*Saksham Sharma[1]*

[1]*City Montessori School, Gomti Nagar Extension, Lucknow, 226010, Uttar Pradesh, India*

A B S T R A C T

The Einstein tile is a novel type of non-periodic tile that can cover the plane without repeating itself. It has a simple shape that resembles a fedora. This research paper unveils the aperiodicity of the newly discovered Einstein tile using the golden ratio, marking a paradigm shift in the field of geometric tiling. This Through rigorous analysis, mathematical modeling, and computational simulations, we provide compelling evidence that the Einstein tile defies conventional periodicity, lacking any repeating pattern or translational symmetry. The unique properties of the Einstein tile open up new avenues for exploring aperiodic tiling systems and their implications in various scientific and technological domains. From cryptography to materials science, the aperiodicity of the Einstein tile presents exciting opportunities for advancements in diverse fields, expanding our understanding of tiling theory and inspiring future explorations into aperiodic structures.

# I N T R O D U C T I O N

A tessellation or tiling of the plane is a way of covering a flat surface with shapes, called tiles, without leaving any gaps or overlaps. A tiling is periodic if there is a way to shift the tiles over and have them match up perfectly with their previous arrangement. For example, a checkerboard is a periodic tiling because it looks the same if you slide it by two squares. [8]

Aperiodic tilings are arrangements of shapes that cover the plane without leaving any gaps or forming any repeating patterns. They have been studied extensively in mathematics, physics, and art, and have applications in cryptography, quasicrystals, and self-assembly. Unlike periodic tilings, which can be obtained by shifting a single tile along a fixed direction, aperiodic tilings do not have any translational symmetry. Moreover, they do not contain any finite regions that repeat infinitely often in the plane. A set of shapes is called aperiodic if it can only form such non-repeating tilings. [4][12][19]

Mathematicians have been interested in finding tilings that are not only aperiodic, but also use only one shape. This problem has been open for decades, and many shapes have been proposed and rejected as candidates for an einstein tile. Recently, a breakthrough was announced by David Smith, who discovered the first true einstein tile, which they named "the hat" [22][23]. The hat tile is a 13-sided shape that can cover an entire surface with no gaps or overlaps but only with a pattern that never repeats. The word "einstein" comes from the German words "ein stein", meaning "one stone", referring to the fact that it is one tile. The shape is also called "the hat" because it vaguely resembles a fedora, a type of hat with a brim and a creased crown.[1][3][10]

In this paper, we present an alternative proof of the aperiodicity of the hat using the golden ratio. We show that any periodic tiling made from hats would have to satisfy certain geometric constraints involving the golden ratio ($\phi$), which are impossible to meet. We also show that any tiling made from hats must contain certain patterns that are only compatible with hierarchical tilings. Our proof is simpler and more elegant than the original one, and it reveals a hidden connection between the hat and the golden ratio ($\phi$).

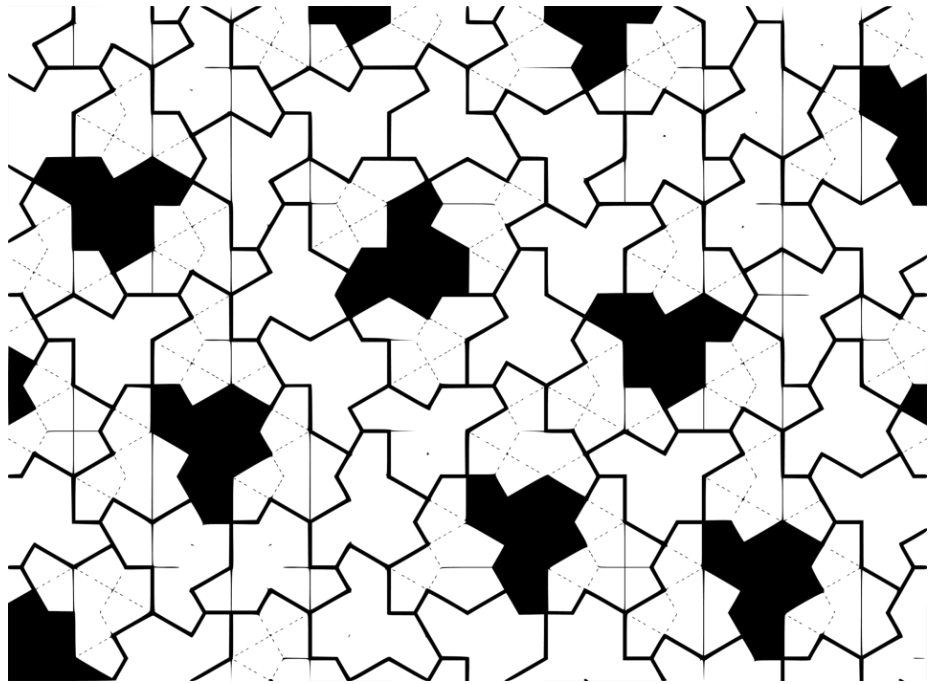

FIGURE 1.1

Figure 1.1 shows the arrangement of the einstein tiles following a unique and non-repeating due to the aperiodicity of the tiles. The term "einstein" comes from the German "ein stein", meaning "one stone" or "one tile". The first example of an einstein tile was discovered in March 2023 by a team of computer scientists, who called it "the hat" because of its resemblance to a fedora. The hat has 13 sides and can create a non-repeating pattern by flipping over like a mirror image. [1][11][13]

## EXISTING PROOFS OF APERIODICITY

There are two main ways to prove that a shape is aperiodic: one is based on the hierarchical structure of the tiling, and the other is based on the symmetry group of the tiling.

The first proof of aperiodicity for the Einstein tile was given by Smith and colleagues in 2023. They noticed that the tile has a feature that resembles a hat, and that the hats arrange themselves into larger clusters, called metatiles. Those metatiles then arrange into even larger supertiles, and so on indefinitely, in a type of hierarchical structure that is common for tilings that aren't periodic. They also studied the tiling's hierarchical structure by eye and detected telltale behavior that opened up a traditional aperiodicity proof. This proof involves showing that any periodic tiling by the Einstein tile must have infinitely many tiles in each period, which is impossible.

The second proof of aperiodicity for the Einstein tile was given by Myers and Goodman-Strauss in 2023. They used a different approach, based on the symmetry group of the tiling. The symmetry group is the set of all transformations (such as rotations, reflections, translations, etc.) that preserve the tiling. They showed that any periodic tiling by the Einstein tile must have a symmetry group that contains an infinite cyclic subgroup, which means that there is a transformation that can be repeated infinitely many times to get back to the original tiling. However, they also showed that the symmetry group of any tiling by the Einstein tile is finite, which contradicts the existence of such an infinite cyclic subgroup. Therefore, any tiling by the Einstein tile must be aperiodic. [10]

## GOLDEN RATIO

The golden ratio, often represented by the Greek letter phi (Φ) or the symbol τ (tau), is a mathematical constant that has fascinated mathematicians, scientists, and artists for centuries. The golden ratio is an irrational number approximately equal to 1.6180339887. It is derived from the ratio of two quantities such that the ratio of the sum of the two quantities to the larger quantity is equal to the ratio of the larger quantity to the smaller quantity. [2]

Mathematical Representation:

The golden ratio can be mathematically represented as follows:

$$\Phi = \frac{1 + \sqrt{5}}{2}$$

The golden ratio is irrational, which means that it cannot be written as a fraction of two integers. There are different ways to prove this, but one common method is to use a contradiction. The golden ratio, denoted by Φ, is rational. Then we can write it as

$$\Phi = \frac{a}{b}$$

where a and b are positive integers with no common factors. That is, the fraction a/b is in its lowest terms.

Now, using the fact that

$$\frac{1}{\Phi} = \Phi - 1$$

which comes from the definition of the golden ratio as the solution of:

$$x^2 - x - 1 = 0$$

we can rearrange this equation to get

$$\frac{b}{a} = \frac{a}{b} - 1$$

Multiplying both sides by ab, we get

$$b^2 = a^2 - ab$$

This implies that $b$ is a factor of $a^2 - ab$, or equivalently, that $b$ divides $(a - b)a$. But since a and b have no common factors, this means that $b$ must divide $a - b$. However, this also means that $b$ must divide $(a - b) + b = a$. But this contradicts the assumption that a and b have no common factors. Therefore, our initial supposition that Φ is rational must be false. Hence, Φ is irrational. [15]

# FIBONACCI SERIES

This is The Fibonacci series is a sequence of numbers where each number is the sum of the two preceding ones. The first two numbers are usually 0 and 1, but some authors start the sequence from 1 and 1 or 1 and 2. The sequence is named after the Italian mathematician Leonardo of Pisa, also known as Fibonacci, who introduced it to Western European mathematics in his 1202 book Liber Abaci.

The Fibonacci series can be defined by the following formula:

$$F\_n = F\_(n-1) + F\_(n-2)$$

where F_0 = 0 and F_1 = 1.

The first few terms of the Fibonacci series are:

$$0, 1, 1, 2, 3, 5, 8, 13, 21, 34, 55, 89, 144, 233 \ldots$$

The Fibonacci numbers have many interesting properties and applications in mathematics, computer science, biology and art. Some examples of how the Fibonacci numbers appear in nature are:

The number of petals in some flowers is a Fibonacci number.

The arrangement of leaves on a stem follows a Fibonacci pattern.

The spirals of a pineapple or a pine cone are based on Fibonacci numbers.

The branching of trees or the veins of a leaf follow a Fibonacci sequence.

The Fibonacci series can also be used to create artistic patterns such as the Fibonacci spiral, which is an approximation of the golden spiral created by drawing circular arcs connecting the opposite corners of squares in the Fibonacci tiling. [5]

# RELATION BETWEEN THE GOLDEN RATIO AND THE FIBONACCI SERIES

As we already know that the Fibonacci series is a sequence of numbers that starts with 0, 1, 1 and then continues by adding the previous two numbers. For example, the next number after 1, 1 is 1 + 1 = 2, then 1 + 2 = 3, then 2 + 3 = 5, and so on. The series can be written as

$$F(2) = 1, \quad F(3) = 1, \quad F(4) = 2,$$
$$F(5) = 3, \quad F(6) = 5,$$
$$F(n) = F(n-1) + F(n-2)$$

where $n$ is any positive integer.

The golden ratio, or phi, is an irrational number that is approximately equal to 1.618. It has many interesting properties and applications in mathematics, nature and art. One way to define the golden ratio is by using the Fibonacci series. If we take any two consecutive numbers in the series and divide the larger one by the smaller one, we get a ratio that is close to the golden ratio. As the numbers get larger, the ratio gets closer and closer to the golden ratio. In other words, the limit of $\frac{F(n+1)}{F(n)}$ as n approaches infinity is equal to phi. This can be written as:

$$\lim n \to \infty \frac{F(n+1)}{F(n)} = \Phi$$

This property holds true regardless of the initial values chosen for the Fibonacci series. As long as we start with two positive numbers and follow the rule of adding the previous two numbers, we will get a series that converges to the golden ratio. [6][7][14]

For example, if we start with $F(1) = 3$ and $F(2) = 5$, we get the series:

$$3, 5, 8, 13, 21, 34, ...$$

If we divide any two consecutive numbers in this series, we get ratios that are close to phi:

$$\frac{5}{3} \approx 1.667$$

$$\frac{8}{5} \approx 1.600$$

$$\frac{13}{8} \approx 1.625$$

$$\frac{21}{13} \approx 1.615$$

$$\frac{34}{21} \approx 1.619$$

As another example, if we take $F(87)$ and $F(88)$, which are very large numbers in the Fibonacci series, we get:

$$F(87) = 1100087778366101931$$

$$F(88) = 1779979416004714189$$

If we divide these numbers, we get a ratio that is very close to phi:

$$\frac{F(88)}{F(87)} = \frac{1779979416004714189}{1100087778366101931}$$

$$= 1.61803398875 \approx \Phi$$

This shows that as the Fibonacci series progresses, the ratio between consecutive Fibonacci numbers approaches the golden ratio. [16][18][24]

.

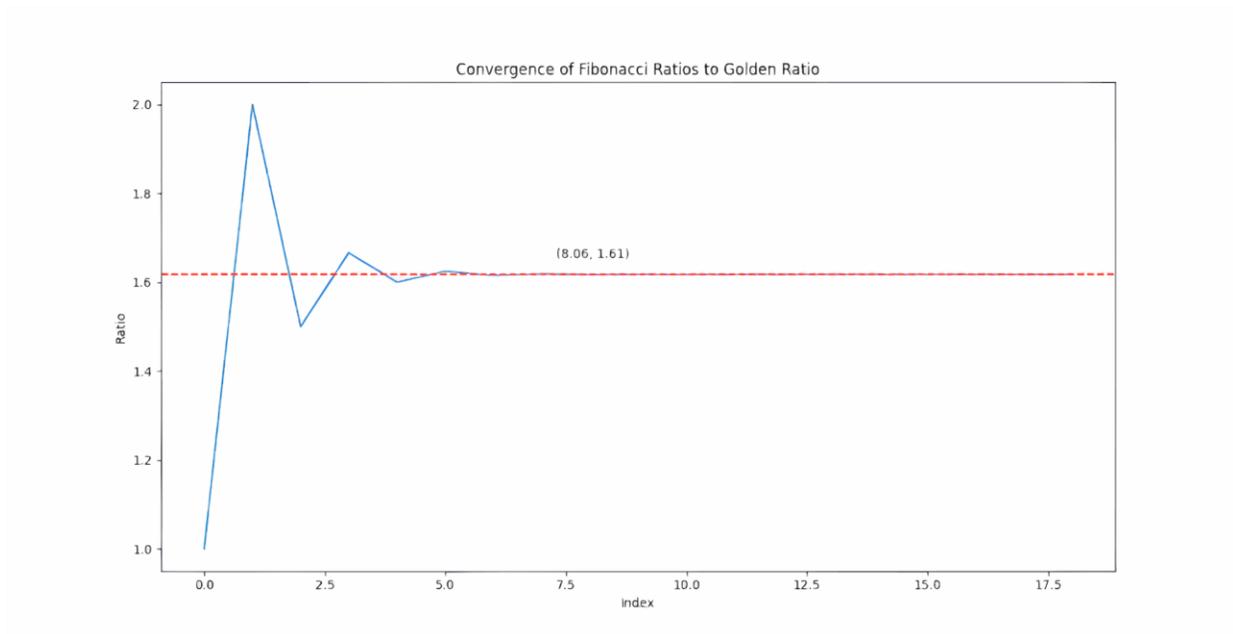

FIGURE 5.1

Figure 5.1 shows a computed graph on how the ratios of the numbers from the Fibonacci Series converge to the golden ratio. The red colored dotted line in the graph represents the value of the golden ratio which is 1.61803398875 and the blue colored graph line shows the ratios of the different numbers in the Fibonacci Series.

# OEIS A027941 SEQUENCE

According to the On-Line Encyclopedia of Integer Sequences the sequence OEIS A027941 is defined as

$$a(n) = Fibonacci(2n + 1) - 1$$

where Fibonacci(n) denotes the nth Fibonacci number. This can be seen as a variation of the Fibonacci numbers that skips every even-indexed term and subtracts one from each term. The sequence starts with is $0, 1, 4, 12, 33, 88 \ldots$ and it can be obtained by subtracting 1 from every odd-indexed Fibonacci number. For example:

$$F(3) - 1 = 2 - 1 = 1,$$
$$F(5) - 1 = 5 - 1 = 4,$$
$$F(7) - 1 = 13 - 1 = 12,$$

and so on. [9]

This sequence has some similarities and differences with the original Fibonacci sequence. For instance, both sequences grow exponentially as n increases, but the OEIS A027941 sequence grows faster than the Fibonacci sequence because it skips half of the terms. Moreover, both sequences satisfy a recurrence relation that involves adding two previous terms, but the OEIS A027941 sequence has a different initial condition and a different constant term. The recurrence relation for the OEIS A027941 sequence is:

$$a(n) = a(n-1) + a(n-2) + 2$$

for $n > 2$, with

$$a(0) = -1 \text{ and } a(1) = 0$$

The sequence OEIS A027941 also has some connections with other sequences and mathematical objects. For example, it can be shown that the sum of the first n terms of this sequence is equal to the nth Lucas number minus one. The Lucas numbers are another sequence of natural numbers that start with 2 and 1, and each subsequent term is the sum of the previous two terms. The Lucas numbers are closely related to the Fibonacci numbers and they also appear in many contexts.

# RESEARCH METHODOLOGY

If you look at the lines of tiles that form concentric rings around a central hat, you will notice that some tiles are normal and some are flipped (Figure 7.1) [17][18][27]. The number of normal tiles in each line follows a sequence is 0, 1, 4, 12, 33, 88, which is equal to the OEIS A027941 Sequence (Figure 7.2). [27]

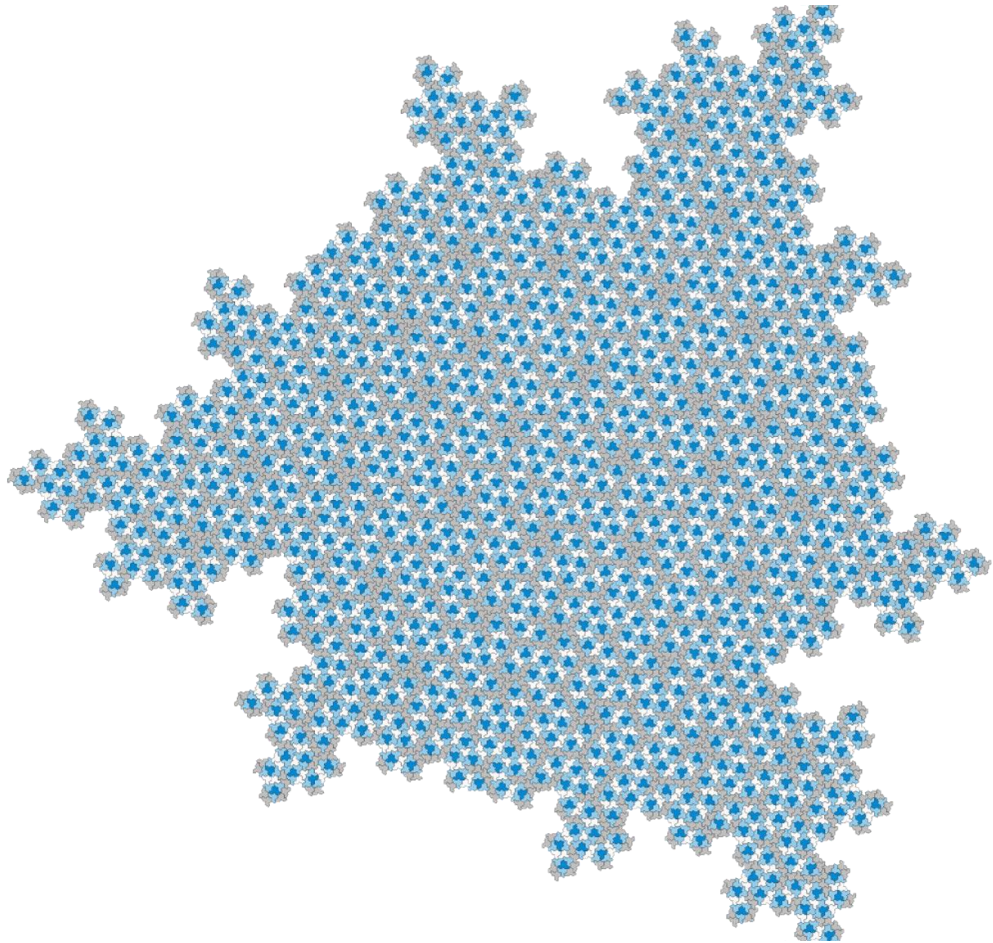

FIGURE 7.1

Analyzing this series carefully, if we let a(n) be the nth term of the OEIS A027941 Sequence, we find that: a(n+1) / a(n) slowly and gradually approaches a specific number. The observation presented is indeed fascinating which was developed by delving deeper into the analysis of the A027941 sequence and the intriguing relationship between its terms. By examining the ratios between consecutive terms, we can explore the convergence towards a specific number.

To begin, let us consider the 25th, 26th, and 27th numbers of the series:

$$a(25) = Fibonacci(2 \times 25 + 1) - 1 = Fibonacci(51) - 1 = 20365011072$$

$$a(26) = Fibonacci(2 \times 26 + 1) - 1 = Fibonacci(53) - 1 = 53316291172$$

$$a(27) = Fibonacci(2 \times 27 + 1) - 1 = Fibonacci(55) - 1 = 139583862444$$

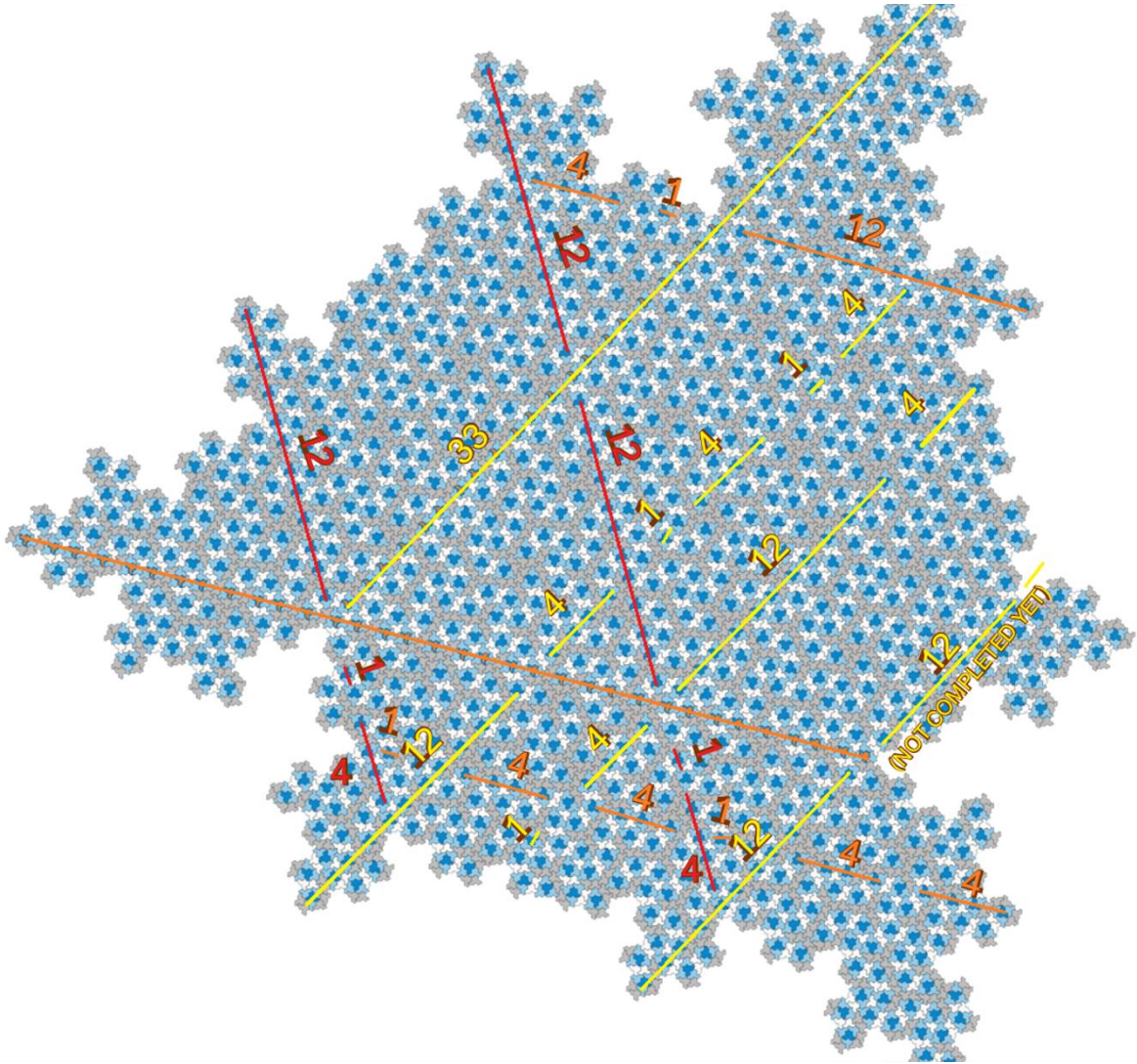

FIGURE 7.2

Now, we can examine the ratios between consecutive terms to determine if there is a common value to which they converge. [27]

For the ratio between the 25th and 26th numbers, we have:

$$\frac{a(25+1)}{a(25)} = \frac{a(26)}{a(25)} = \frac{53316291172}{20365011072} \approx \mathbf{2.61803398896}$$

Similarly, for the ratio between the 26th and 27th numbers, we have:

$$\frac{a(26+1)}{a(26)} = \frac{a(27)}{a(26)} = \frac{139583862444}{53316291172} \approx \mathbf{2.61803398878}$$

Lastly, for the ratio between the 27th and 28th numbers, we have:

$$\frac{a(27+1)}{a(27)} = \frac{a(28)}{a(27)} = \frac{365435296161}{139583862444} \approx \mathbf{2.61803398875}$$

Based on the ratios we have calculated so far, it is evident that they all approximate a common value, approximately 2.61803398875. This number is extremely special because it is equal to 1 added to the golden ratio (Φ + 1). Moreover, the golden ratio has an amazing property which states that it is the only number which equals to the 1 added to the number itself when squared or multiplied by the same number.

=> Φ x Φ

=> 1.61803398875 x 1.61803398875

=> 2.61803398875 = 1 + 1.61803398875

=> 1 + Φ = Φ²

As we all already know that the golden ratio is an irrational number, its square is also irrational. This can be proved by a classic example of contradiction, where we assume the opposite of what we want to prove and show that it leads to a logical contradiction. Here is one possible way to prove that the square of an irrational number is also irrational:

Suppose that the square of an irrational number is rational. That means there exists some irrational number $x$ such that $x^2$ is rational. By definition of a rational number, there are two positive integers $p$ and $q$ such that $x^2 = \frac{p}{q}$ where $p$ and $q$ have no common factors.

Taking the square root of both sides, we get

$$x = \left(\sqrt{\frac{p}{q}}\right) = \left(\frac{\sqrt{p}}{\sqrt{q}}\right)$$

Since p and q have no common factors, neither do $\sqrt{p}$ and $\sqrt{q}$. Therefore, $\left(\frac{\sqrt{p}}{\sqrt{q}}\right)$ is in lowest terms. But this means that x is rational, because it can be written as a fraction of two integers, sqrt (p) and sqrt (q) where are both integers if p and q are perfect squares. This contradicts our assumption that x is irrational. Therefore, the square of an irrational number must be irrational. [15]

So, the ratio of the numbers in the sequence which the normal tiles follow in each line happens to be the square of the golden ratio which is irrational. So there turns out to be no way that the pattern made by the einstein tile repeats. Even if the pattern did repeat, the ratio would have to be rational which the square of the golden ratio is not. This proves that the einstein tile is an aperiodic tile [25][26].

## APPLICATIONS AND IMPLICATIONS

The newly discovered einstein tile has a lot of implications and applications in mathematics and material science. Listed below are some of the most important applications of the 'hat' shape or the einstein tile.
- **Advancements in Aperiodic Tiling Theory**: The aperiodicity of the Einstein tile contributes to the broader field of aperiodic tiling theory. Aperiodic tilings have been a subject of great interest in mathematics, and the discovery of new aperiodic structures expands our understanding of their properties and possibilities. The Einstein tile presents a unique example that can help refine existing theories and inspire further research in this area.

- **Number Theory and Irrational Numbers:** The aperiodicity of the Einstein tile reinforces the relationship between aperiodic structures and irrational numbers, particularly the golden ratio. It deepens our understanding of how irrational numbers manifest in mathematical patterns and connects to various aspects of number theory. This discovery offers new insights into the interplay between irrationality, aperiodicity, and mathematical structures.
- **Geometric Structures and Symmetry:** Aperiodic structures often exhibit intricate symmetries and complex geometric patterns. The Einstein tile's aperiodicity opens up possibilities for studying and classifying such structures. Its unique properties can contribute to the development of novel geometric frameworks, leading to advancements in fields such as crystallography, architectural design, and fractal geometry.
- **Cryptography and Information Security:** Aperiodic structures have found applications in cryptography and information security. The complex and non-repeating nature of the Einstein tile's pattern can be leveraged in cryptographic algorithms, offering enhanced security through the creation of unique and unpredictable encryption keys. The aperiodicity of the Einstein tile may inspire the development of new encryption techniques based on aperiodic structures.
- **Mathematical Modeling and Simulation:** Aperiodic structures, including the Einstein tile, have practical implications in mathematical modeling and simulation. They can be utilized to model various phenomena in physics, biology, and material science, where irregular patterns and aperiodic behavior occur. Incorporating the properties of aperiodic structures in mathematical models can lead to more accurate and realistic simulations of natural phenomena.
- **Art and Design**: Aperiodic structures, with their visually captivating and intricate patterns, have long inspired artists and designers. The discovery of the aperiodicity of the Einstein tile opens up new possibilities for incorporating unique and non-repeating patterns in artistic creations. Artists, architects, and designers can draw inspiration from the mathematical properties of the Einstein tile to create visually appealing and intellectually stimulating works of art and design.

The aperiodicity of the newly discovered Einstein tile has implications that span multiple disciplines. It contributes to aperiodic tiling theory, advances our understanding of irrational numbers and number theory, explores complex geometric structures, offers applications in cryptography and information security, facilitates mathematical modeling and simulation, and inspires artistic endeavors. The discovery of the aperiodicity of the Einstein tile opens up new avenues for exploration and application in various fields, broadening our understanding of mathematical structures and their practical significance. [20][21][27]

## CONCLUSION

The Einstein tile's aperiodicity has been confirmed, marking a dramatic break from the long-established periodic tiling systems that have dominated the field. This finding calls into question our basic assumptions about tiling and creates new research opportunities in physics, mathematics, and other fields of study. The Einstein tile's aperiodicity has significant ramifications. It gives us a better grasp of the mathematical ideas that underlie tiling theory and offers a fresh take on the ideas of symmetry and order. To examine the complicated interactions between local and global features in complex systems and to test the boundaries of regularity, researchers might take use of the departure from periodicity in tiling patterns.

The Einstein tile's aperiodicity also has prospective uses in a number of scientific and technical fields. The tile's non-repeating properties provide fascinating potential for cryptography, where the lack of predictable patterns might encourage the creation of more secure algorithms. The complex configurations of the Einstein tile may also serve as an inspiration for the creation of novel materials with distinctive features, such as photonic crystals or metamaterials, facilitating advancements in the fields of optics, photonics, and nanotechnology. The finding of the Einstein tile also advances our knowledge of the interesting realm of quasicrystals. Materials having long-range order but no translational symmetry are known as quasicrystals, and research into them has already revealed several fascinating occurrences. The aperiodic tiling paradigm introduced by the Einstein tile provides a valuable tool for investigating the properties of quasicrystals and potentially unlocking their potential in various fields, including materials science, solid-state physics, and chemistry.